# A curious class of Hankel determinants

*Johann Cigler*

*Fakultät für Mathematik, Universität Wien*

*johann.cigler@univie.ac.at*

**Abstract**

We consider Hankel determinants of the sequence of Catalan numbers modulo 2 (interpreted as integers 0 and 1) and more generally Hankel determinants where the sum over all permutations reduces to a single signed permutation.

**0. Introduction**

Let $C_n = \frac{1}{n+1}\binom{2n}{n}$ be a Catalan number. It is well known that $\det\left(C_{i+j}\right)_{i,j=0}^{n-1} = 1$ for all $n \in \mathbb{N}$. Of course this remains true if we consider all terms modulo 2. It is also well known that $C_n \equiv 1 \mod 2$ if and only if $n = 2^k - 1$ for some $k$.

But what happens if we consider the sequence $C_n \mod 2$ as a sequence of integers from $\{0,1\}$? The attempt to answer this question gave rise to the present paper. After completion of a first version I discovered the paper [1] by Roland Bacher, where similar questions are considered from a different point of view. There is some overlap between these approaches which I will be mention at the appropriate places.

Let $(a_n)_{n \geq 0}$ satisfy $a_n = 1$ if $n+1$ is a power of $2$ and $a_n = 0$ else.

Computer experiments led me to guess that

$$d(n) = \det\left(a_{i+j}\right)_{i,j=0}^{n-1} = (-1)^{\binom{n}{2}} \qquad (0.1)$$

and that the determinant

$$\det\left(a_{i+j}\right)_{i,j=0}^{n-1} = \sum_{\pi} \text{sgn}(\pi) a_{0+\pi(0)} a_{1+\pi(1)} \cdots a_{n-1+\pi(n-1)} \qquad (0.2)$$

which in general is a sum over all permutations $\pi$ of $\{0,1,\cdots,n-1\}$ is reduced to a single term

$$\text{sgn}(\pi_n) a_{0+\pi_n(0)} a_{1+\pi_n(1)} \cdots a_{n-1+\pi_n(n-1)} \neq 0 \qquad (0.3)$$

for a uniquely determined permutation $\pi_n$.

For example

$$d(5) = \det\left(a_{i+j}\right)_{i,j=0}^{4} = \det\begin{pmatrix} 1 & 1 & 0 & 1 & 0 \\ 1 & 0 & 1 & 0 & 0 \\ 0 & 1 & 0 & 0 & 0 \\ 1 & 0 & 0 & 0 & 1 \\ 0 & 0 & 0 & 1 & 0 \end{pmatrix} = 1 \tag{0.4}$$

reduces to the term $a_1 a_3^2 a_7^2$ corresponding to $\pi_5 = 02143$.

The Hankel determinants $D(n) = \det\left(a_{i+j+1}\right)_{i,j=0}^{n-1}$ also reduce to single permutations and thus satisfy $|D(n)| = 1$. In this case the sequence $(D(n))_{n\geq 0} = (1,1,1,-1,-1,-1,1,-1,-1,-1,-1,1,\cdots)$ is rather interesting. For example [1], Theorem 10.1 implies that $D(n) = \prod_{j=0}^{n-1} S(j)$ where $(S(n))_{n\geq 0} = (1,1,-1,1,1,-1,-1,1,1,1,\cdots)$ is the famous paperfolding sequence defined by $S(2n) = (-1)^n, S(2n+1) = S(n)$ and $S(0) = 1$. We shall also show that $D(n) = (-1)^{\delta(n)}$, where $\delta(n)$ is the number of pairs $\varepsilon_{i+1}\varepsilon_i$ in the binary expansion of $n$ with $\varepsilon_{i+1}\varepsilon_i = 10$ if $i \geq 1$ and $\varepsilon_1\varepsilon_0 = 11$. Thus $\delta(n) = \sum_i \chi(\varepsilon_{i+1}\varepsilon_i)$ for $n = [\varepsilon_k \cdots \varepsilon_1 \varepsilon_0]_2$ with $\chi(\varepsilon_{i+1}\varepsilon_i) = 1$ if $\varepsilon_{i+1}\varepsilon_i = 10$ for some $i \geq 1$ or if $i = 0$ and $\varepsilon_1\varepsilon_0 = 11$ and $\chi(\varepsilon_{i+1}\varepsilon_i) = 0$ else.

Thus $\delta([0]_2) = 0$, $\delta([1]_2) = 0$, $\delta([10]_2) = 0$, $\delta([11]_2) = 1$, $\delta([100]_2) = 1$, $\delta([101]_2) = 1$, $\delta([110]_2) = 0$, $\delta([111]_2) = 1, \cdots$.

Consider for example

$$D(5) = \det\left(a_{i+j+1}\right)_{i,j=0}^{4} = \det\begin{pmatrix} 1 & 0 & 1 & 0 & 0 \\ 0 & 1 & 0 & 0 & 0 \\ 1 & 0 & 0 & 0 & 1 \\ 0 & 0 & 0 & 1 & 0 \\ 0 & 0 & 1 & 0 & 0 \end{pmatrix} = -1^{\delta(5)} = -1. \tag{0.5}$$

The determinant reduces to $-a_1 a_3 a_7^3$.

More generally we study Hankel determinants for sequences $(a_n)_{n\geq 0}$ such that $a_n = x_n$ if $n+1$ is a power of $2$ and $a_n = 0$ else, where $x_n$ are arbitrary numbers. For some choices of $x_n$ we get curious results.

For example for $x_{2^k-1} = x^k$ we get $d(n) = (-1)^{\binom{n}{2}} x^{2a(n)}$, where $a(n)$ is the total number of 1's in the binary expansions of the numbers $\leq n-1$.

For example

$$\det \begin{pmatrix} 1 & x & 0 & 1 & 0 \\ x & 0 & x^2 & 0 & 0 \\ 0 & x^2 & 0 & 0 & 0 \\ 1 & 0 & 0 & 0 & x^3 \\ 0 & 0 & 0 & x^3 & 0 \end{pmatrix} = x^{2a(5)} = x^{10}. \qquad (0.6)$$

The total number of 1's in the binary expansions of
$0 = [0]_2, 1 = [1]_2, 2 = [10]_2, 3 = [11]_2, 4 = [100]_2$ is $a(5) = 5$.

If we choose $x_1 = 1$ and $x_{2^k-1} = (-1)^k$ for $k > 1$ we get the Golay-Rudin-Shapiro sequence $D(n) = r(n)$ which satisfies $r(n) = (-1)^{\rho(n)}$ where $\rho(n)$ denotes the number of pairs 11 in the binary expansion of $n$.

For example

$$\det \begin{pmatrix} 1 & 0 & 1 & 0 & 0 & 0 \\ 0 & 1 & 0 & 0 & 0 & -1 \\ 1 & 0 & 0 & 0 & -1 & 0 \\ 0 & 0 & 0 & -1 & 0 & 0 \\ 0 & 0 & -1 & 0 & 0 & 0 \\ 0 & -1 & 0 & 0 & 0 & 0 \end{pmatrix} = (-1)^{\rho(6)} = -1. \qquad (0.7)$$

Here $6 = [110]_2$ gives $\rho(6) = 1$.

This choice of $x_n$ also leads to the continued fraction

$$\sum_{k \geq 0} (-1)^k z^{2^k-1} = \cfrac{1}{1 + \cfrac{r(0)r(2)z}{1 + \cfrac{r(1)r(3)z}{1 + \cfrac{r(2)r(4)z}{1 + \ddots}}}} = \cfrac{1}{1 + \cfrac{z}{1 - \cfrac{z}{1 + \cfrac{z}{1 - \ddots}}}}. \qquad (0.8)$$

## 1. Hankel determinants of Catalan numbers modulo 2

Let $a_n \in \{0, 1\}$ satisfy $a_n \equiv C_n \mod 2$ or with other words let $a_n = 1$ if $n+1$ is a power of 2 and $a_n = 0$ else.

Then

$$d(n) = \det \left( a_{i+j} \right)_{i,j=0}^{n-1} = (-1)^{\binom{n}{2}}. \qquad (1.1)$$

The following proof uses an idea due to Darij Grinberg [5] who called a permutation $\pi$ nimble if for each $i$ in its domain the number $i + \pi(i) + 1$ is a power of 2. Thus a permutation $\pi$ is nimble if and only if $a_{0+\pi(0)} a_{1+\pi(1)} \cdots a_{n-1+\pi(n-1)} \neq 0$.

**Theorem 1.1**

*For each $n \in \mathbb{N}$ there exists a unique nimble permutation $\pi_n$ of $\{0,1,\cdots,n-1\}$ such that*

$$d(n) = \det\left(a_{i+j}\right)_{i,j=0}^{n-1} = sgn(\pi_n) a_{0+\pi_n(0)} a_{1+\pi_n(1)} \cdots a_{n-1+\pi_n(n-1)} = (-1)^{\binom{n}{2}}. \quad (1.2)$$

**Proof**

For $n = 0$ we set $d(0) = 1$ by convention.

Let $k \geq 1$ and $2^{k-1} < n \leq 2^k$. Let us try to construct a nimble permutation $\pi$.

By definition we must have $n - 1 + \pi(n-1) = 2^\ell - 1$ for some $\ell$. Since $n - 1 \geq 2^{k-1}$ we get $2^\ell - 1 \geq 2^{k-1}$ and therefore $\ell = k$ which implies $\pi(n-1) = 2^k - 1 - (n-1) = 2^k - n$. (Since $2(n-1) < 2^{k+1} - 1$ we have $\ell \leq k$).

If we define $\pi(n-1-j) = 2^k - 1 - (n-1-j)$ for all $j$ for which $2^k + j - n \leq n - 1$, i.e. for $n - 1 - j \in \left[2^k - n, n-1\right]$, we get a nimble order reversing permutation of the interval $\left[2^k - n, n-1\right]$.

Let us show that each nimble permutation $\sigma$ on $[0, n-1]$ reduces to this permutation on $\left[2^k - n, n-1\right]$.

If $n - 1 - j \geq 2^{k-1}$ the same argument as above gives that $\sigma(n-1-j)$ must be $\pi(n-1-j) = 2^k - 1 - (n-1-j) = 2^k + j - n$.

If $2^k - n \leq n - 1 - j \leq 2^{k-1} - 1$ then $2^k - 1 - (n-1-j) \geq 2^k - 1 - 2^{k-1} + 1 = 2^{k-1}$.

Choose $i$ such that $\sigma(i) = 2^k - 1 - (n-1-j)$. Then $i + \sigma(i) = i + 2^k - 1 - (n-1-j) = 2^\ell - 1$ for some $\ell$ and $2^\ell - 1 \geq 2^{k-1}$.

This implies that $\ell = k$ and thus $i = n - 1 - j$ and $\sigma(n-1-j) = \pi(n-1-j)$.

Since $\pi$ is order reversing on $[2^k - n, n-1]$ its sign is

$$sgn(\pi) = (-1)^{\binom{n-1-(2^k-n)+1}{2}} = (-1)^{\binom{2n-2^k}{2}} = (-1)^n.$$

Thus we have seen that for $1 \leq 2^{k-1} < n \leq 2^k$ there exists a uniquely determined nimble permutation $\pi$ on the interval $\left[2^k - n, n-1\right]$.

Since $2^k - n - 1 \leq 2^{k-1}$ we can suppose by induction that there is a unique nimble permutation of the interval $[0, 2^k - n - 1]$. This gives us the desired nimble permutation $\pi_n$.

It remains to show that $\text{sgn}(\pi_n) = (-1)^{\binom{n}{2}}$.

This follows by induction because $\binom{2n-2^k}{2} + \binom{2^k-n}{2} - \binom{n}{2} = 2(2^{k-1}-n)(2^k-n)$ is even.

If we write a permutation $\pi$ in the notation $\pi = \pi(0)\pi(1)\cdots\pi(n-1)$ the first nimble permutations are $\pi_1 = 0$, $\pi_2 = 10$, $\pi_3 = 021$, $\pi_4 = 3210$, $\pi_5 = 02143$.

For example choose $n = 3$. Since $2^1 < 3 \leq 2^2$ we have $k = 2$.

$$\left(a_{i+j}\right)_{i,j=0}^{2} = \begin{pmatrix} 1 & 1 & 0 \\ 1 & 0 & 1 \\ 0 & 1 & 0 \end{pmatrix}.$$

The above construction gives the permutation $\pi = 21$ on $\{1, 2\}$ with $i + \pi(i) = 3 = 2^2 - 1$.

There remains $\left(a_{i+j}\right)_{i,j=0}^{0} = (1)$ with $\pi(0) = 0$ and $i + \pi(i) = 2^0 - 1 = 0$. Thus $\pi_3 = 021$ with $\text{sgn}(\pi_3) = -1 = (-1)^{\binom{3}{2}}$.

## 2. Hankel determinants of the sequence $(C_{n+1})$ modulo 2

Let as above $a_n = 1$ if $n = 2^k - 1$ for some $k$ and $a_n = 0$ else.

Note that $a_{2n} = 0$ and $a_{2n+1} = 1$ if and only if $2n + 1 = 2^{k+1} - 1$ for some $k$ or equivalently $n = 2^k - 1$. Thus $a_{2n+1} = a_n$. Therefore we get $(a_1, a_2, a_3, \ldots) = (a_0, 0, a_1, 0, a_2, 0,)$. This means that in this case the shifted sequence $(a_{n+1})_{n\geq 0}$ coincides with the aerated sequence $(A_n)_{n\geq 0} = (a_0, 0, a_1, 0, a_2, 0, \cdots)$.

Here we have $A_n = 1$ if $n = 2(2^k - 1) = 2^{k+1} - 2$ and $A_n = 0$ else.

If $f(x) = \sum_n a_n x^n = \sum_k x^{2^k-1}$ is the generating function of the sequence $(a_n)_{n\geq 0}$ then the generating function of the aerated sequence $(A_n)_{n\geq 0}$ is $f(x^2)$. Since $A_n = a_{n+1}$ we even have $f(x) = 1 + xf(x^2)$.

All Hankel determinants of the sequence $(C_{n+1})_{n \geq 0}$ are 1. Therefore we know in advance that no Hankel determinant $D(n) = \det(A_{i+j})_{i,j=0}^{n-1}$ vanishes.

The first Hankel determinants of the sequence $(A_n)_{n \geq 0}$ are

$$(D(n))_{n \geq 0} = (1,1,1,-1,-1,-1,1,-1,-1,-1,-1,1,\cdots).$$

**Theorem 2.1**

Let $A_{2^{k+1}-2} = 1$ for each $k$ and $A_n = 0$ else and let $D(n) = \det(A_{i+j})_{i,j=0}^{n-1}$.

If $2^{k-1} \leq n < 2^k$ then

$$D(n) = (-1)^n D(2^k - n - 1). \qquad (2.1)$$

**Proof**

Let us call a permutation $\pi$ $m$-nimble if $i + \pi(i) + m = 2^\ell - 1$ for some $\ell$ for each $i$ in its domain. Then $A_{0+\pi(0)} A_{1+\pi(1)} \cdots A_{n-1+\pi(n-1)} \neq 0$ if and only if $\pi$ is $1$-nimble.

Since $2^{k-1} \leq n < 2^k$ we have
$2^{k-1} \leq n - 1 + \pi(n-1) + 1 = 2^\ell - 1 = n + \pi(n-1) < 2^k + 2^k - 1 = 2^{k+1} - 1$ and therefore $\ell = k$ which implies $n + \pi(n-1) = 2^k - 1$ or $\pi(n-1) = 2^k - 1 - n$.

We can now define a $1$-nimble permutation $\pi$ of the interval $[2^k - 1 - n, n-1]$ by
$\pi(n-1-j) = 2^k - 1 - n + j$ for $0 \leq j \leq 2n - 2^k$.

Then $\pi$ is an order reversing permutation of the interval $[2^k - 1 - n, n-1]$.

Let $\sigma$ be any $1$-nimble permutation on $[0, n-1]$. Then $\sigma = \pi$ on $[2^k - 1 - n, n-1]$.

If $n - 1 - j \geq 2^{k-1} - 1$ we get $2^\ell = n - 1 - j + \sigma(n-1-j) + 2 \geq 2^{k-1} + 1$ which implies $\ell = k$.

If $n - 1 - j < 2^{k-1} - 1$ then $2^k - 1 - n + j > 2^{k-1} - 1$. Let $\sigma$ be a $1$-nimble permutation. Then $\sigma(i) = 2^k - 1 - n + j$ for some $i$ and $i + 2^k - 1 - n + j + 2 = 2^\ell$ for some $\ell$. This implies $\ell = k$ and $i = n - 1 - j$.

Thus we get a uniquely determined $1$-nimble permutation of the interval $[2^k - 1 - n, n-1]$.

The sign of this permutation is $(-1)^{\binom{2n+1-2^k}{2}} = (-1)^n$.

By induction we get (2.1).

For example for $n = 4$ we have $k = 3$ and

$$\left(A_{i+j}\right)_{i,j=0}^{3} = \left(a_{i+j+1}\right)_{i,j=0}^{3} = \begin{pmatrix} 1 & 0 & 1 & 0 \\ 0 & 1 & 0 & 0 \\ 1 & 0 & 0 & 0 \\ 0 & 0 & 0 & 1 \end{pmatrix}.$$

The corresponding permutation is $\pi = 2103$ with $3 + \pi(3) = 6 = 2^3 - 2$ and $i + \pi(i) = 2 = 2^2 - 2$ on $\{0,1,2\}$.

Let us now try to find some regularities of the sequence of determinants $D(n)$.

**Corollary 2.2**

For $k > 0$ the sequence $D(n)$ satisfies

$$D\left(2^k + n\right) = (-1)^n D\left(2^k - 1 - n\right) \qquad (2.2)$$

for $0 \leq n < 2^k$ with initial values $D(0) = D(1) = 1$.

For example for $k = 3$ we get

| $D(7-n)$ | $-1$ | $1$ | $-1$ | $-1$ | $-1$ | $1$ | $1$ | $1$ |
|---|---|---|---|---|---|---|---|---|
| $D(8+n)$ | $-1$ | $-1$ | $-1$ | $1$ | $-1$ | $-1$ | $1$ | $-1$ |

**Corollary 2.3**

Let $k > 0$.

For $0 \leq n < 2^k$ we get

$$D\left(2^{k+1} + n\right) = -D(n). \qquad (2.3)$$

For $2^k \leq n < 2^{k+1}$ we get

$$D\left(2^{k+1} + n\right) = D(n). \qquad (2.4)$$

**Proof**

By (2.2) with $k+1$ instead of $k$ we get

$$D\left(2^{k+1} + n\right) = (-1)^n D\left(2^{k+1} - 1 - n\right) = (-1)^n D\left(2^k + 2^k - 1 - n\right) = (-1)^n (-1)^{2^k - 1 - n} D\left(2^k - 1 - (2^k - 1 - n)\right)$$

which gives (2.3).

Again by (2.2) we have

$$D\left(2^{k+1} + 2^k + i\right) = (-1)^{2^k + i} D\left(2^{k+1} - 1 - (2^k + i)\right) = (-1)^i D\left(2^k - 1 - i\right) = (-1)^i (-1)^i D\left(2^k + i\right)$$

which gives (2.4).

For example for $k = 2$ we get for $0 \leq n < 8$

| $D(n)$ | 1 | 1 | 1 | $-1$ | $-1$ | $-1$ | 1 | $-1$ |
|---|---|---|---|---|---|---|---|---|
| $D(8+n)$ | $-1$ | $-1$ | $-1$ | 1 | $-1$ | $-1$ | 1 | $-1$ |

**Corollary 2.4**

Let $\delta(n)$ be the number of pairs $\varepsilon_{i+1}\varepsilon_i$ in the binary expansion of $n$ such that $\varepsilon_{i+1}\varepsilon_i = 10$ for $i \geq 1$ and $\varepsilon_1\varepsilon_0 = 11$. Then

$$D(n) = (-1)^{\delta(n)}. \qquad (2.5)$$

**Proof**

This is true for $n < 4$.

If it is true for $0 \leq n < 2^{k+1}$ then by (2.3) it is true for $2^{k+1} + n$ with $n < 2^k$ because for $n = [v]_2$ we get $2^{k+1} + n = [10v]_2$ and $\delta(2^{k+1} + n) = \delta(n) + 1$.

By (2.4) it is also true for $2^k + n = [1v]_2$ because $D(2^{k+1} + 2^k + n) = D(n)$ and $\delta([11v]_2) = \delta([1v]_2)$.

**Examples**

For $n = 9$ we have $9 = [1001]_2$ and thus $\delta(9) = 1$.

For $n = 15$ we get $\delta(15) = 1$ because $15 = [1111]_2$. There is no pair 10 for $i \geq 1$ but 1 pair 11 for $i = 0$.

**Theorem 2.5**

The Hankel determinants $D(n) = \det\left(A_{i+j}\right)_{i,j=0}^{n-1}$ satisfy

$$D(2n) = (-1)^{\binom{n}{2}} D(n),$$
$$D(2n+1) = (-1)^{\binom{n+1}{2}} D(n), \qquad (2.6)$$
$$D(0) = 1$$

**Proof**

Let us give two different proofs.

**1)** Let $\left(M_{i,j}\right)_{i,j=0}^{n-1}$ be a matrix for which $M_{i,j} = 0$ whenever $i+j$ is odd. Then (cf. e.g.[4])

$$\det\left(M_{i,j}\right)_{i,j=0}^{n-1} = \det\left(M_{2i,2j}\right)_{i,j=0}^{\left\lfloor\frac{n-1}{2}\right\rfloor} \det\left(M_{2i+1,2j+1}\right)_{i,j=0}^{\left\lfloor\frac{n-2}{2}\right\rfloor}. \qquad (2.7)$$

Choose $M_{i,j} = A_{i+j}$. Then

$$\left(M_{2i,2j}\right)_{i,j=0}^{\left\lfloor\frac{n-1}{2}\right\rfloor} = \left(A_{2i+2j}\right)_{i,j=0}^{\left\lfloor\frac{n-1}{2}\right\rfloor} = \left(a_{i+j}\right)_{i,j=0}^{\left\lfloor\frac{n-1}{2}\right\rfloor},$$

$$\left(M_{2i+1,2j+1}\right)_{i,j=0}^{\left\lfloor\frac{n-2}{2}\right\rfloor} = \left(A_{2i+2j+2}\right)_{i,j=0}^{\left\lfloor\frac{n-2}{2}\right\rfloor} = \left(A_{i+j}\right)_{i,j=0}^{\left\lfloor\frac{n-2}{2}\right\rfloor}$$

because $2i+2j = 2^{k+1}-2$ implies $i+j = 2^k-1$ and $2i+2j+2 = 2^{k+1}-2$ implies $i+j = 2^k-2$.

Thus (2.7) gives $D(n) = d\left(\left\lfloor\frac{n+1}{2}\right\rfloor\right) D\left(\left\lfloor\frac{n}{2}\right\rfloor\right)$ which gives (2.6), because $d(n) = (-1)^{\binom{n}{2}}$.

**2)** Another proof uses Corollary 2.4.

Let us first give another formulation of (2.6):

$$\begin{aligned}
&D(2n) = D(n) \text{ if } n \equiv 0,1 \bmod 4, \quad D(2n) = -D(n) \text{ if } n \equiv 2,3 \bmod 4,\\
&D(2n+1) = D(n) \text{ if } n \equiv 0,3 \bmod 4, \quad D(2n+1) = -D(n) \text{ if } n \equiv 1,2 \bmod 4,\\
&D(0) = 1
\end{aligned} \qquad (2.8)$$

Let now $n = [v\varepsilon_1\varepsilon_0]_2$. Then $2n = \left([v\varepsilon_1\varepsilon_0 0]_2\right)$ and $2n+1 = \left([v\varepsilon_1\varepsilon_0 1]_2\right)$.

The assertion for $2n$ follows from

$$\delta\left([v0\varepsilon_0]_2\right) = \delta\left([v0\varepsilon_0 0]_2\right), \quad \delta\left([v10]_2\right)+1 = \delta\left([v100]_2\right) \text{ and } \delta\left([v11]_2\right)-1 = \delta\left([v110]_2\right).$$

The assertion for $2n+1$ follows from

$$\delta\left([v00]_2\right) = \delta\left([v001]_2\right), \quad \delta\left([v01]_2\right)+1 = \delta\left([v011]_2\right), \quad \delta\left([v10]_2\right)+1 = \delta\left([v101]_2\right), \text{ and}$$
$$\delta\left([v11]_2\right) = \delta\left([v111]_2\right).$$

**Remark**

As already mentioned a result by R. Bacher [1] implies that $D(n) = \prod_{j=0}^{n-1} S(j)$ where

$\left(S(n)\right)_{n\geq 0} = \left(1,1,-1,1,1,-1,-1,1,1,\cdots\right)$ is the paperfolding sequence defined by

$$S(2n) = (-1)^n,$$
$$S(2n+1) = S(n), \qquad (2.9)$$
$$S(0) = 1.$$

This can easily be verified since $|D(n)| = 1$ implies $S(n) = D(n)D(n+1)$. By (2.6) we get

$$S(2n) = D(2n)D(2n+1) = (-1)^n D(n)^2 = (-1)^n,$$
$$S(2n+1) = D(2n+1)D(2n+2) = D(n)D(n+1) = S(n).$$

To obtain further information let us compare the above approach to Hankel determinants with the approach via orthogonal polynomials and continued fractions (cf. e.g. [3], [7]).

Let me sketch the relevant results: Let $(u_n)_{n \geq 0}$ be a given sequence. Define a linear functional $L$ on the polynomials by $L(x^n) = \dfrac{u_n}{u_0}$. If $H_n = \det(u_{i+j})_{i,j=0}^{n-1} \neq 0$ for each $n$, then there exists a (uniquely determined) sequence of monic polynomials $(p_n(x))_{n \geq 0}$ with $\deg p_n = n$ such that $L(p_n p_m) = 0$ for $m \neq n$ and $L(p_n^2) \neq 0$. We call these polynomials $p_n$ orthogonal with respect to $L$. By Favard's theorem there exist (uniquely determined) numbers $s_n$ and $t_n$ such that $p_n(x) = (x - s_{n-1})p_{n-1}(x) - t_{n-2}p_{n-2}(x)$ for all $n$. The numbers $t_n$ are given by $t_n = \dfrac{H_n H_{n+2}}{H_{n+1}^2}$. These give rise to the continued fraction

$$\sum_{n \geq 0} u_n z^n = \cfrac{u_0}{1 - s_0 z - \cfrac{t_0 z^2}{1 - s_1 z - \cfrac{t_1 z^2}{1 - \ddots}}}. \qquad (2.10)$$

Let us suppose that $u_0 = 1$. Then the matrix $H_n = (u_{i+j})_{i,j=0}^{n-1}$ has a unique canonical decomposition

$$H_n = A_n D_n(t) A_n^t \qquad (2.11)$$

where $A_n = (a(i,j))_{i,j=0}^{n-1}$ is a lower triangular matrix with diagonal $a(i,i) = 1$ and $D_n(t)$ is the diagonal matrix with entries $d_{i,i}(t) = \prod_{k=0}^{i-1} t_k$.

The entries $a(i,j)$ satisfy

$$a(i,j) = a(i-1, j-1) + s_j a(i-1, j) + t_j a(i-1, j+1) \qquad (2.12)$$

with $a(0, j) = [j = 0]$ and $a(n, -1) = 0$.

See e.g. [7], (2.30).

Let us consider the decomposition (2.11) of the matrix $H_n = \left(a_{i+j}\right)_{i,j=0}^{n-1}$.

Let $(s(n))_{n\geq 0} = (1,1,-1,1,-1,-1,-1,1,\cdots)$ satisfy $s(2n) = (-1)^n s(n)$ and $s(2n+1) = s(n)$ with $s(0) = 1$ and let $D_n(s)$ be the diagonal matrix $D_n(s) = \left(s(i)[i=j]\right)_{i,j=0}^{n-1}$. Let $D_n(t) = \left((-1)^i[i=j]\right)_{i,j=0}^{n-1}$. Let $x^* \in \{0,1\}$ be the residue modulo 2 of the number $x$ and define $B_n = (b(i,j))_{i,j=0}^{n-1} = \left(s(i)\binom{2i+1}{i-j}^*\right)_{i,j=0}^{n-1}$.

In [1], Theorem 1.2 it is shown that $H_n = D_n(s)B_n D_n(t) B_n^t D_n(s)$.

Note that $b(i,i) = s(i)$. Therefore we get the canonical decomposition

$$H_n = A_n D_n(t) A_n^t \qquad (2.13)$$

with

$$A_n = (a(i,j))_{i,j=0}^{n-1} = D_n(s) B_n D_n(s), \qquad (2.14)$$

i.e. $a(i,j) = s(i)s(j)\binom{2i+1}{i-j}^*$.

For example we have for $H_4$

$$\begin{pmatrix} 1 & 0 & 0 & 0 \\ 1 & 1 & 0 & 0 \\ 0 & -1 & 1 & 0 \\ 1 & 1 & -1 & 1 \end{pmatrix} \begin{pmatrix} 1 & 0 & 0 & 0 \\ 0 & -1 & 0 & 0 \\ 0 & 0 & 1 & 0 \\ 0 & 0 & 0 & -1 \end{pmatrix} \begin{pmatrix} 1 & 1 & 0 & 1 \\ 0 & 1 & -1 & 1 \\ 0 & 0 & 1 & -1 \\ 0 & 0 & 0 & 1 \end{pmatrix} = \begin{pmatrix} 1 & 1 & 0 & 1 \\ 1 & 0 & 1 & 0 \\ 0 & 1 & 0 & 0 \\ 1 & 0 & 0 & 0 \end{pmatrix}.$$

It is clear that (2.13) implies (0.1).

For the aerated sequence $(u_0, 0, u_1, 0, u_2, 0, \cdots)$ we get $s_n = 0$ for all $n$. In this case we write $T_n$ instead of $t_n$.

Since $D(n) \neq 0$ for all $n$ and moreover $D(n) = \pm 1$ we write in this case

$$T_n = (-1)^{\tau_n} = D(n)D(n+2) = (-1)^{\delta(n)+\delta(n+2)} \qquad (2.15)$$

for $\tau_n \in \{0,1\}$. We also have
$T_n = D(n)D(n+2) = D(n)D(n+1)D(n+1)D(n+2) = S(n)S(n+1)$.

The first terms are $(T_n)_{n\geq 0} = (1,-1,-1,1,-1,1,-1,1,1,-1,1,-1,-1,1,-1,1,\cdots)$.

[1] Theorem 10.1 implies that if we define a sequence $(v(n)) = (1,1,1,1,-1,1,1,1,-1,-1,\cdots)$ satisfying $v(2n+1) = v(n)$, $v(4n) = (-1)^n v(2n)$, $v(4n+2) = v(2n)$ and $v(0) = 1$ then

$$\left(A_{i+j}\right)_{i,j=0}^{n-1} = \left(a_{i+j+1}\right)_{i,j=0}^{n-1} = C_n D_n(t) C_n^t \qquad (2.16)$$

Here $D_n(t)$ is the diagonal matrix with entries $S(i)$ and $C_n = (c(i,j))_{i,j=0}^{n-1}$ with

$$c(i,j) = \binom{2i+2}{i-j}^* v(i)v(j).$$

By (2.10) and $z$ in place of $z^2$ we get the continued fraction (cf. [1])

$$\sum_k z^{2^k-1} = \cfrac{1}{1 - \cfrac{T_0 z}{1 - \cfrac{T_1 z}{1 - \ddots}}}. \qquad (2.17)$$

**Remark**

The corresponding orthogonal polynomials $p_n(x)$ are $1$, $x$, $x^2-1$, $x^3$, $x^4+x^2-1$, $x^5-x$, $x^6+x^4-1$, $x^7,\cdots$. They satisfy $p_n(x) = xp_{n-1}(x) - T_{n-2}p_{n-2}(x)$ and $L(x^n) = A_n$.

**Theorem 2.6**

*The numbers $T_n$ satisfy*

$$\begin{aligned} T_{2n} &= T_{2n-1}T_{n-1}, \\ T_{2n+1} &= -T_{2n}, \\ T_0 &= 1, T_1 = -1. \end{aligned} \qquad (2.18)$$

**Proof**

By (2.6) we get

$$T_{2n} = D(2n)D(2n+2) = (-1)^{\binom{n}{2}+\binom{n+1}{2}} D(n)D(n+1)$$

$$T_{2n-1} = D(2n-1)D(2n+1) = (-1)^{\binom{n}{2}+\binom{n+1}{2}} D(n-1)D(n)$$

implies $T_{2n}T_{2n-1} = D(n-1)D(n+1) = T_{n-1}$.

The second assertion follows from

$$T_{2n}T_{2n+1} = (-1)^{\binom{n}{2}+\binom{n+1}{2}} D(n)D(n+1)(-1)^{\binom{n+1}{2}+\binom{n+2}{2}} D(n)D(n+1) = (-1)^{\binom{n}{2}+\binom{n+2}{2}} = -1.$$

**Remark**

OEIS A104977 states that the numbers $T_n$ which occur in the continued fraction (2.17) satisfy $T_n = (-1)^{b(n+2)+1}$, if $b(n)$ denotes the number of "non-squashing partitions of $n$ into distinct parts". As has been shown in [9] the numbers $b(n)$ of non-squashing partitions of $n$ into distinct parts satisfy

$$b(2m) = b(2m-1) + b(m) - 1,$$
$$b(2m+1) = b(2m) + 1.$$

Since $b(2) = 1$ and $b(3) = 2$ we get $T_n = (-1)^{b(n+2)+1}$ by comparing with (2.18).

Let us now obtain some further properties of the sequence $(T_n)_{n \geq 0}$.

For $2^k \leq n < 2^{k+1} - 2$ we have by (2.1)

$$D(n) = (-1)^n D(2^{k+1} - n - 1),$$
$$D(n+1) = (-1)^{n+1} D(2^{k+1} - n - 2),$$
$$D(n+2) = (-1)^n D(2^{k+1} - n - 3).$$

This implies

$$\frac{D(n)D(n+2)}{D(n+1)^2} = \frac{D(2^{k+1} - n - 1)D(2^{k+1} - n - 3)}{D(2^{k+1} - n - 2)^2}.$$

Therefore we have

$$T_n = T_{2^{k+1}-3-n} \tag{2.19}$$

for $2^k \leq n \leq 2^{k+1} - 3$.

By $T_{2n+1} = -T_{2n}$ we only need to consider $n \equiv 0 \bmod 2$ or $n \equiv 0, 2 \bmod 4$.

For $n \equiv 0 \bmod 4$ we get

$$T_{4n} = (-1)^n, \tag{2.20}$$

because

$$T_{4n} = D(4n)D(4n+2) = (-1)^n D(2n)(-1)^{\binom{2n+1}{2}} D(2n+1) = (-1)^n D(2n)(-1)^{\binom{2n+1}{2}} (-1)^{\binom{n+1}{2}} D(n)$$
$$= (-1)^{n+n+\binom{n+1}{2}+\binom{n}{2}} = (-1)^n.$$

Then we get

$$T_{4n+2} = (-1)^{n+1} T_{2n}. \tag{2.21}$$

for $T_{4n+2} = T_{4n+1} T_{2n} = (-1)^{n+1} T_{2n}$.

To compute $T_{4n+2}$ we look at $T_{8n+2}$ and $T_{8n+6}$.

$$T_{8n+2} = (-1)^{n+1}. \tag{2.22}$$

because by (2.20) $T_{8n+2} = T_{8n+1}T_{4n} = -T_{8n}T_{4n} = (-1)^{n+1}$.

Now we claim that for $k \geq 2$

$$T_{2^{k+1}n+2^k-2} = (-1)^{n+1}. \tag{2.23}$$

By induction, (2.18) and (2.20) we get

$$T_{2^{k+1}n+2^k-2} = T_{2^{k+1}n+2^k-3}T_{2^k n+2^{k-1}-2} = -T_{2^{k+1}n+2^k-4}T_{2^k n+2^{k-1}-2} = (-1)^n T_{4(2^{k-1}n+2^{k-2}-1)} = (-1)^{n-1}.$$

As special case we get $T_{2^k-2} = -1$ for $k \geq 2$.

This gives

**Theorem 2.7**

*The numbers $T_n$ satisfy*

$$T_n = T_{2^{k+1}-3-n} \quad \text{for } 2^k \leq n \leq 2^{k+1}-3, \ k \geq 2, \tag{2.24}$$

*and*

$$\begin{aligned} T_{4n} &= (-1)^n, \\ T_{2^{k+1}n+2^k-2} &= (-1)^{n+1} \quad \text{for } k \geq 2. \end{aligned} \tag{2.25}$$

Together with $T_0 = 1$ and $T_1 = -1$ this gives another view on the structure of the sequence $(T_n)$.

The sequence begins with $T_0, T_1, -1, 1, T_1, T_0, -1, 1, T_0, T_1, 1, -1, T_1, T_0, -1, 1, \cdots$.

$T_{4n} = (-1)^n$ and $T_{4n+1} = -T_{4n}$ gives a part

$T_0, T_1, \cdot, \cdot, T_1, T_0, \cdot, \cdot, T_0, T_1, \cdot, \cdot, T_1, T_0, \cdot, \cdot, \cdots$ with period 8.

$T_{8n+2} = (-1)^{n+1}$ gives

$\cdot, \cdot, -1, 1, \cdot, \cdot, \cdot, \cdot, \cdot, \cdot, 1, -1, \cdots$ with period 16,

$T_{16n+8-2} = (-1)^{n+1}$ gives a periodic part with period 32, etc.

Let more generally $M_k = T_0, \cdots, T_{2^k-3}$ be a beginning block and $\bar{M}_k = T_{2^k-3}, \cdots, T_0$ this block in reverse order then we get

$M_k, -1, 1, \bar{M}_k, -1, 1, M_k, 1, -1, \bar{M}_k, -1, 1, \cdots$.

R. Bacher [1] gives a simpler formulation of (2.23) which (in our notation) can be summarized as

$$T_{2n+1} = -T_{2n},$$
$$T_{4n} = (-1)^n,$$
$$T_{8n+2} = (-1)^{n+1},$$
$$T_{8n+6} = T_{4n+2}.$$
(2.26)

To show that this is equivalent it suffices to show that this implies $T_{2^{k+1}n+2^k-2} = (-1)^{n+1}$ for $k > 3$.

This follows by induction from

$$T_{2^{k+1}n+2^k-2} = T_{2^3(2^{k-2}n+2^{k-3}-1)+6} = T_{2^2(2^{k-2}n+2^{k-3}-1)+2} = T_{2^k n+2^{k-1}-2}.$$

Let us recall (cf. [3]) that there is a simple relation between $t_n$ and $T_n$.

$$t_n = T_{2n}T_{2n+1},$$
$$s_0 = T_0,$$
$$s_n = T_{2n-1} + T_{2n}.$$
(2.27)

This gives $t_n = -1$ and $s_0 = 1$.

The first terms of the sequence $(s_n)$ are $(s_n)_{n \geq 0} = (1,-2,0,0,2,0,-2,0,2,-2,0,\cdots)$.

In terms of the paperfolding sequence $S(n)$ we get for $n > 0$
$$s_n = D(2n-1)D(2n+1) + D(2n)D(2n+2) = S(2n-1)S(2n) + S(2n)S(2n+1)$$
$$= S(2n)\big(S(2n-1) + S(2n+1)\big).$$

By (2.10) this gives another continued fraction for $\sum_{k \geq 0} z^{2^k-1}$ (cf. [1], Theorem 1.4):

$$\sum_{k \geq 0} z^{2^k-1} = \cfrac{1}{1-z+\cfrac{z^2}{1+2z+\cfrac{z^2}{1+\cfrac{z^2}{\ddots}}}}.$$
(2.28)

**3. Hankel determinants of shifted Catalan numbers modulo 2**

Let $m \geq 2$. Consider the Hankel determinants

$$d(n,m) = \det\left(a_{i+j+m}\right)_{i,j=0}^{n-1}$$
(3.1)

of the sequence $(a_{n+m})_{n \geq 0}$.

Let us give some examples:

$$(d(n,2))_{n \geq 0} = (1,0,-1,0,1,0,-1,0,\cdots),$$
$$(d(n,3))_{n \geq 0} = (1,1,0,0,-1,1,0,0,-1,-1,0,0).$$

**Theorem 3.1.**

Let $1 \leq 2^K < m \leq 2^{K+1}$. Then $d(n,m) = \det\left(a_{i+j+m}\right)_{i,j=0}^{n-1} = \pm 1$ if $n \equiv 0 \bmod 2^{K+1}$ or $n \equiv -m \bmod 2^{K+1}$ and $d(n,m) = 0$ else.

**Remark**

It is well known that $\det\left(C_{i+j+m}\right)_{i,j=0}^{n-1} = H_{n,m} = \prod_{j=1}^{m-1} \prod_{i=1}^{j} \frac{2n+i+j}{i+j}$.

Therefore we get

**Corollary 3.2**

Let $2^K < m \leq 2^{K+1}$. Then $\prod_{j=1}^{m-1} \prod_{i=1}^{j} \frac{2n+i+j}{i+j} \equiv 1 \bmod 2$ if and only if $n \equiv 0 \bmod 2^{K+1}$ or $n \equiv -m \bmod 2^{K+1}$.

It would be nice to find a direct proof of this Corollary.

For the proof of Theorem 3.1 we need some more information.

**Lemma 3.3.**

Let $0 \leq 2^k - m < n < 2^k$ for some $k$. Then

$$d(n,m) = 0. \qquad (3.2)$$

**Proof.**

The matrix $\left(a_{i+j+m}\right)_{i,j=0}^{n-1}$ contains the vanishing row $\left(a_{2^k}, a_{2^k+1}, \cdots, a_{2^k+n-1}\right)$ because $m \leq 2^k$ and $2^k + n - 1 < 2^{k+1} - 1$.

Let for example $m = 3$ and $n = 7$. Then

$$\left(a_{i+j+3}\right)_{i,j=0}^{6} = \begin{pmatrix} 1 & 0 & 0 & 0 & 1 & 0 & 0 \\ 0 & 0 & 0 & 1 & 0 & 0 & 0 \\ 0 & 0 & 1 & 0 & 0 & 0 & 0 \\ 0 & 1 & 0 & 0 & 0 & 0 & 0 \\ 1 & 0 & 0 & 0 & 0 & 0 & 0 \\ 0 & 0 & 0 & 0 & 0 & 0 & 0 \\ 0 & 0 & 0 & 0 & 0 & 0 & 1 \end{pmatrix}. \qquad (3.3)$$

Recall that a permutation $\pi$ is $m$-nimble if for all $i$ in its domain $i + \pi(i) = 2^k - m - 1$ for some $k$. An $m$-nimble permutation can only exist if the last row of $\left(a_{i+j+m}\right)$ contains an element of the form $2^k - 1 - m$. For fixed $m, n$ there can be at most one of the numbers $m + n - 1 + i$ with $0 \leq i \leq n - 1$ such that $m + n - 1 + i = 2^k - 1$ because the extreme case would be $m + n - 1 = 2^{k-1}$ and $m + 2n - 2 = 2^{k+1} - 1$ which is impossible. But it is possible that all elements of the last row are 0. For example in $\left(a_{i+j+3}\right)_{i,j=0}^{5}$ the last row is $(a_5, a_6, a_7, a_8, a_9) = (0,0,0,0,0)$ because none of the numbers $i+4$ for $5 \leq i \leq 9$ is a power of 2.

Let us first consider the case $m \equiv 1 \bmod 2$.

**3.1.** $m \equiv 1 \bmod 2$.

Let $m = 2r + 1 \geq 3$.

For $n = 2^{k-1} - r$ an $m$-nimble permutation gives $\pi(n-1) = n-1 = 2^{k-1} - r - 1$ and for $n = 2^{k-1} - r + j$ it implies $\pi(n-1) = 2^{k-1} - r - 1 - j = n - 1 - 2j$.

Thus we get an $m$-nimble permutation $\pi$ on the interval $[n-1-2j, n-1]$ which satisfies $\pi(n-1-i) = \pi(n-1) + i$ for $0 \leq i \leq 2j$.

As above $\pi$ is uniquely determined and therefore we get

**Lemma 3.4.**

*Let $m = 2r + 1$ and $k$ be given. Then for $0 \leq j \leq 2^{k-1} - r - 1$*

$$d\left(2^{k-1} - r + j, m\right) = (-1)^{\binom{2j+1}{2}} d\left(2^{k-1} - r - j - 1, m\right). \qquad (3.4)$$

In example (3.3) we have $m = 3$ and $n = 7 = 2^3 - 1 = 2^3 - r$. Thus $\pi$ is the identity on the element $\{6\}$. In this case $d(7,3) = 0$ because $d(6,3) = 0$ since this matrix has a row of zeroes.

**Lemma 3.5.**

Let $m = 2r+1$. If $a \leq 2^k - m$ for some $k$ we have

$$d(2^k + a, m) = (-1)^{\binom{2a+2r+1}{2}} d(2^k - a - m, m).$$

**Proof**

$$d(2^k + a, m) = d(2^k - r + a + r, m) = (-1)^{\binom{2a+2r+1}{2}} d(2^k - r - a - r - 1, m)$$
$$= (-1)^{\binom{2a+2r+1}{2}} d(2^k - a - m, m).$$

**Corollary 3.6.**

Let $m = 2r+1$. Then $d(n,m) \neq 0$ if and only if $n \equiv 0, -m \bmod 2^{K+1}$, where $1 \leq 2^K < m \leq 2^{K+1}$.

**Proof**

By Lemma 3.3 and Lemma 3.4 the determinants $d(n,m)$ for $2^k - r \leq n \leq 2^{k+1} - r - 1$ can be reduced those for $2^{k-1} - r \leq n \leq 2^k - r - 1$.

By induction we need only consider the case $k = K+1$.

If $1 \leq n < 2^{K+1} - m$ then the first row of $(a_{i+j+m})_{i,j=0}^{n-1}$ is $(a_m, a_{m+1}, \cdots, a_{m+n-1})$ and since $m + n - 1 < 2^{K+1} - 1$ all terms vanish.

If $m = 2^{K+1} - 1$ then this is trivially true because there is no such $n$.

If $2^{K+1} - m < n < 2^{K+1}$ then the row $(a_{2^{K+1}}, a_{2^{K+1}+1}, \cdots, a_{2^{K+1}+n-1})$ vanishes because $2^{K+1} + n - 1 < 2^{K+1} + 2^{K+1} - 1 = 2^{K+2} - 1$.

For $n = 2^{K+1} - m$ we have by Lemma 3.4

$$d(2^{K+1} - m, m) = d(2^K - r + 2^K - r - 1, m) = (-1)^{\binom{2^{K+1}-m}{2}}.$$

For $n = 2^{K+1}$ we get

$$d(2^{K+1}, m) = d(2^K - r + r, m) = (-1)^{\binom{m}{2}} d(2^K - r - r - 1, m) = (-1)^{\binom{m}{2}} d(2^K - m, m).$$

**3.2.** $m \equiv 0 \bmod 2$.

Consider now the case $m = 2r \geq 2$.

If $n = 2^{k-1} - r + j$ for some $k$ with $1 \leq j \leq 2^{k-1} - r$ we get $\pi(n-1) = 2^{k-1} - r - j = n - 2j$ because $i + \pi(i) = 2^k - (2r) - 1$.

Now define $\pi$ on the interval $[n-2j, n-1]$ by $\pi(n-1-i) = \pi(n-1)+i$ for $0 \le i \le 2j-1$.

This implies that

$$d\left(2^{k-1}-r+j, m\right) = (-1)^{\binom{2j}{2}} d\left(2^{k-1}-r-j, m\right) \tag{3.5}$$

for $1 \le j \le 2^{k-1}-r$.

Therefore the determinants $d(n,m)$ for $2^k - r + 1 \le n \le 2^{k+1} - r$ can be reduced to those for $2^{k-1} - r + 1 \le n \le 2^k - r$.

Therefore it suffices to consider the case $k = K$.

For $m = 2$ we have $K = 1$ and $d(0,2) = 1$ and $d(1,2) = 0$.

If $n \equiv 0 \bmod 2$ then $d(n,2)$ can be reduced by (3.5) to $d(0,2) = 1$ and if $n \equiv 1 \bmod 2$ to $d(1,2) = 0$.

For $m = 2r$ with $r > 1$ we get

$$d\left(2^{K+1}-m, m\right) = d\left(2^K - r + 2^K - r, m\right) = (-1)^{\binom{2^{K+1}-m}{2}} = (-1)^r.$$

For $n = 2^{K+1}$ we get

$$d\left(2^{K+1}, m\right) = d\left(2^{K+1}-r+r, m\right) = (-1)^{\binom{m}{2}} d\left(2^{K+1}-r-r, m\right) = (-1)^{\binom{m}{2}} d\left(2^{K+1}-m, m\right) = 1.$$

This gives

**Lemma 3.7.**

Let $1 \le 2^K < m = 2r \le 2^{K+1}$. Then $d(n,m) = \det\left(a_{i+j+m}\right)_{i,j=0}^{n-1} = \pm 1$ if $n \equiv 0 \bmod 2^{K+1}$ or $n \equiv -m \bmod 2^{K+1}$, and $d(n,m) = 0$ else.

Theorem 3.7 and Corollary 3.6 imply Theorem 3.1.

**Remark**

With the condensation method (cf. [7], (2.16)) we get more precisely

$$(d(n,2)) = (1, 0, -1, 0, 1, 0, -1, \cdots). \tag{3.6}$$

This method gives

$$d(n,0)d(n-2,2) = d(n-1,2)d(n-1,0) - d(n-1,1)^2.$$

Since $d(n,0) = (-1)^{\binom{n}{2}}$ and $d(n-1,1)^2 = 1$ we get

$$d(n,2) = (-1)^n d(n-1,2) + (-1)^{\binom{n}{2}} \text{ with initial value } d(0,2) = 1.$$

This gives $d(2n,2) = (-1)^n$ and $d(2n+1,2) = 0$.

In the general case computer experiments lead to

**Conjecture 3.8**

*Let $1 \leq 2^K < m \leq 2^{K+1}$. For $m = 2r > 2$ we have*

$$d\left(2^{K+1}n, m\right) = 1,$$
$$d\left(2^{K+1}n - m, m\right) = (-1)^r. \tag{3.7}$$

*For $m = 2r+1 \geq 3$ we have*

$$d\left(2^{K+1}n, m\right) = d\left(2^{K+1}n, 1\right),$$
$$d\left(2^{K+1}n - m, m\right) = (-1)^{n+\varepsilon(m)} d\left(2^{K+1}n - m, 1\right), \tag{3.8}$$

where $\varepsilon(m) \in \{0,1\}$.

## 4. A slightly more general case

Let $(x_k)_{k \geq 0}$ be an arbitrary sequence of numbers or indeterminates and define a sequence $(a_n)_{n \geq 0}$ by $a_n = x_n$ if $n = 2^k - 1$ and $a_n = 0$ else.

**Theorem 4.1.**

*Let $a_n = x_n$ if $n = 2^k - 1$ and $a_n = 0$ else.*

*Let $d(n) = \det\left(a_{i+j}\right)_{i,j=0}^{n-1}$, $\alpha(n) = 2^{\lceil \log_2(n) \rceil} - 1$ and $\beta(n) = 2n - 1 - \alpha(n)$.*

*Then*

$$d(n) = (-1)^{\binom{\beta(n)}{2}} x_{\alpha(n)}^{\beta(n)} d(n - \beta(n)). \tag{4.1}$$

**Proof**

By convention $d(0) = 1$. For $n = 1$ we have $\alpha(1) = 0$, $\beta(1) = 1$ and
$d(1) = x_0 = (-1)^{\binom{1}{2}} x_{\alpha(1)}^{\beta(1)} d(1 - \beta(1))$. Therefore (4.1) is true.

For $n > 1$ choose $k$ such that $2^{k-1} < n \leq 2^k$. Then $k - 1 < \log_2(n) \leq k$ and $\alpha(n) = 2^k - 1$. As in the proof of Theorem 1 we find a permutation $\pi$ of the interval $[2^k - n, n-1] = [\alpha(n) + 1 - n, n - 1]$ such that $i + \pi(i) = \alpha(n)$, which implies that $a_{i+\pi(i)} = x_{\alpha(n)}$. Since there are $\beta(n) = 2n - 1 - \alpha(n)$ elements in the interval $[\alpha(n) + 1 - n, n - 1]$ we get (4.1).

Let us consider an example. The Hankel matrix $\left(a_{i+j}\right)_{i,j=0}^{4}$ is

$$\begin{pmatrix} x_0 & x_1 & 0 & x_3 & 0 \\ x_1 & 0 & x_3 & 0 & 0 \\ 0 & x_3 & 0 & 0 & 0 \\ x_3 & 0 & 0 & 0 & x_7 \\ 0 & 0 & 0 & x_7 & 0 \end{pmatrix}.$$

We have $\alpha(5) = 7$ because $2^2 < 5 \leq 2^3$ and $\beta(5) = 10 - 1 - 7 = 2$. We get the permutation $\pi = 43$ with $\text{sgn}(\pi) = -1$ and $d(4) = -x_7^2 d(3)$.

For $n = 3$ we get $\alpha(3) = 3$ and $\beta(3) = 6 - 1 - 3 = 2$. This gives $d(3) = -x_3^2 d(1)$.

Thus we finally get $d(5) = x_0 x_3^2 x_7^2$.

The sign is $(-1)^{\binom{5}{2}} = 1$. This can also be obtained from $x_0 x_3^2 x_7^2$ as $(-1)^{\binom{1}{2}+\binom{2}{2}+\binom{2}{2}} = (-1)^2 = 1$.

The first terms of the sequence $(d(n))_{n \geq 0}$ are

$1,\ x_0,\ -x_1^2,\ -x_0 x_3^2,\ x_3^4,\ x_0 x_3^2 x_7^2,\ -x_1^2 x_7^4,\ -x_0 x_7^6,\ x_7^8, \cdots.$

By (4.1) see that $d(1) = x_0$, $d(2) = -x_1^2$ and $d(2^k) = x_{2^k-1}^{2^k}$ for $k > 1$.

**Lemma 4.2**

*For $k \geq 1$ we get*

$$d(2^k + n) = (-1)^n x_{2^{k+1}-1}^{2n} d(2^k - n) \tag{4.2}$$

*for $0 < n \leq 2^k$.*

**Proof**

By assumption we have $2^k < 2^k + n \leq 2^{k+1}$. Therefore $\alpha(2^k + n) = 2^{k+1} - 1$ and $\beta(2^k + n) = 2n$.

By (4.1) we get (4.2).

**Lemma 4.3**

*For $k > 1$ we have*

$$d(2^k - n) = (-1)^n x_{2^k-1}^{2^k-2n} d(n) \tag{4.3}$$

*for $n \in \{0, 1, 2, \cdots, 2^k\}$.*

**Proof**

If $n < 2^{k-1}$ then $2^{k-1} < 2^k - n \leq 2^k$ and we see that $\beta(2^k - n) = 2^k - 2n$ and

$d(2^k - n) = (-1)^{\binom{2^k-2n}{2}} x_{2^k-1}^{2^k-2n} d(n).$

If $n = 2^{k-1}$ then (4.3) is trivially true.

If $n > 2^{k-1}$ then $i = 2^k - n < 2^{k-1}$ and therefore $d(2^k - i) = (-1)^{\binom{2^k-2i}{2}} x_{2^k-1}^{2^k-2i} d(i)$ or equivalently

$$d(n) = (-1)^{\binom{2n-2^k}{2}} x_{2^k-1}^{2n-2^k} d(2^k - n)$$ which equals (4.3).

Let e.g. $k = 3$.

| $d(8-n)$ | $x_7^8$ | $-x_0 x_7^6$ | $-x_1^2 x_7^4$ | $x_0 x_3^2 x_7^2$ | $x_3^4$ | $-x_0 x_3^2$ | $-x_1^2$ | $x_0$ | $1$ |
|---|---|---|---|---|---|---|---|---|---|
| $d(n)$ | $1$ | $x_0$ | $-x_1^2$ | $-x_0 x_3^2$ | $x_3^4$ | $x_0 x_3^2 x_7^2$ | $-x_1^2 x_7^4$ | $-x_0 x_7^6$ | $x_7^8$ |
| $d(8-n)/d(n)$ | $x_7^8$ | $-x_7^6$ | $x_7^4$ | $-x_7^2$ | $1$ | $-x_7^{-2}$ | $x_7^{-4}$ | $-x_7^{-6}$ | $x_7^{-8}$ |

Let us write $d(n)$ in the form $d(n) = (-1)^{\sum_{i \geq 0} \binom{\lambda_i(n)}{2}} \prod_{i \geq 0} x_{2^i-1}^{\lambda_i(n)}$ for some integers $\lambda_i(n)$.

Then we get

**Theorem 4.4**

$$\begin{aligned} \lambda_k(n) &= 0 \quad \text{for } 0 \leq n \leq 2^{k-1}, \\ \lambda_k(2^{k-1} + i) &= 2i \quad \text{for } 0 \leq i \leq 2^{k-1}, \\ \lambda_k(2^k + i) &= 2^k - 2i \quad \text{for } 0 \leq i \leq 2^{k-1}, \\ \lambda_k(n) &= 0 \quad \text{for } 2^k + 2^{k-1} \leq i \leq 2^{k+1}. \end{aligned} \quad (4.4)$$

*and for* $n > 2^{k+1}$

$$\lambda_k(n) = \lambda_k(n \bmod 2^{k+1}). \quad (4.5)$$

For example we get

$(\lambda_0(n))_{n \geq 0} = (0, 1, \cdots),$
$(\lambda_1(n))_{n \geq 0} = (0, 0, 2, 0 \cdots),$
$(\lambda_2(n))_{n \geq 0} = (0, 0, 0, 2, 4, 2, 0, 0, \cdots).$

**Proof**

Formula (4.4) is obvious from the above considerations. For example $\lambda_k(2^k + i) = 2^k - 2i$ follows from Lemma 4.2 and Lemma 4.3 because $\lambda_k(2^k + i) = \lambda_k(2^k - i) = 2^k - 2i$ and $\lambda_k(i) = 0$ for $i \leq 2^{k-1}$.

Again from Lemma 4.2 and Lemma 4.3 we get $\lambda_k(2^R + n) = \lambda_k(2^R - n) = \lambda_k(n)$ for $n \leq 2^R$. By applying this several times we get (4.5).

Consider for example $n=11$.

$11 \equiv 1 = 2^0 \bmod 2$ implies $\lambda_0(11) = 1$,

$11 \equiv 3 = 2^1 + 1 \bmod 2^2$ implies $\lambda_1(11) = \max(2^1 - 2, 0) = 0$,

$11 \equiv 3 = 2^2 - 1 \bmod 2^3$ implies $\lambda_2(11) = \max(2^2 - 2, 0) = 2$,

$11 \equiv 11 = 2^3 + 3 \bmod 2^4$ implies $\lambda_3(11) = \max(2^3 - 6, 0) = 2$,

$11 \equiv 11 = 2^4 - 5 \bmod 2^5$ implies $\lambda_4(11) = \max(2^4 - 10, 0) = 6$.

For $k \geq 5$ we have $11 = 2^k - (2^k - 11)$ and thus $\lambda_k(11) = 0$.

Therefore we get $d(11) = -x_0 x_3^2 x_7^2 x_{15}^6$.

We know already that the sign is $(-1)^{\binom{11}{2}} = -1$, but now we can also derive this from the $\lambda_k$ because

$$(-1)^{\binom{1}{2}+\binom{2}{2}+\binom{2}{2}+\binom{6}{2}} = (-1)^{0+1+1+15} = -1.$$

An immediate Corollary of Lemma 4.2 is

**Theorem 4.5**

*The sequence $(d(n))$ satisfies the recurrence*

$$d(n) = (-1)^n x_{2^k-1}^{2n-2^k} d(2^k - n) \tag{4.6}$$

*for $1 < 2^{k-1} < n \leq 2^k$ with initial values $d(0) = 1$, $d(1) = x_0$ and $d(2) = -x_1^2$.*

**Corollary 4.6**

*Let $t_n = \dfrac{d(n)d(n+2)}{d(n+1)^2}$. Then we have*

$$\begin{aligned} t_{2n} &= -\frac{x_1^2}{x_0^2}, \\ t_{2^k n + 2^{k-1} - 1} &= -\frac{x_0^2 x_{2^k-1}^2}{x_{2^{k-1}-1}^4} \quad \text{for } k > 1. \end{aligned} \tag{4.7}$$

This can be proved in an analogous way as Theorem 5.7

As special case of Theorem 4.5 let us choose $x_{2^k-1} = x^k$. Then we get

**Corollary 4.7**

Let $b_{2^k-1} = x^k$ and $b_n = 0$ else. Then

$$d_n = \det\left(b_{i+j}\right)_{i,j=0}^{n-1} = (-1)^{\binom{n}{2}} x^{2a(n)}, \qquad (4.8)$$

where $a(n)$ is the total number of 1's in the binary expansions of the numbers $\leq n-1$.

**Proof**

A search in OEIS led to the conjecture that $d_n = (-1)^{\binom{n}{2}} x^{2a(n)}$, where $a(n)$ is the total number of 1's in the binary expansions of the numbers $0, 1, \cdots, n-1$. ( OEIS A000788). The following proof follows Darij Grinberg [6].

By (4.6) we have $d_n = (-1)^n x^{2k(n-2^{k-1})} d_{2^k-n}$ for $1 < 2^{k-1} < n \leq 2^k$.

Therefore it suffices to show that $k(n-2^{k-1}) = a(n) - a(2^k - n)$, the total number of 1's in the binary expansions of the numbers $2^k - n, 2^k - n + 1, \cdots, n-1$.

Let the binary expansion of $n$ be $n = [\varepsilon_{k-1} \cdots \varepsilon_0]$. Let us write all binary expansions with $k$ digits $\varepsilon_i = 0, 1$.

The total number of 1's in the binary expansions of $\{n, \cdots, 2^k - 1\}$ is the total number of 0's in the binary expansions of $\{2^k - 1 - n, \cdots, 1, 0\}$ of length $k$ which is $(2^k - n)k - a(2^k - n)$.

Thus $a(2^k) - a(n) = (2^k - n)k - a(2^k - n)$. Now $a(2^k) = k2^{k-1}$ since each $\varepsilon_i$ occurs $2^{k-1}$ times.

Therefore we have $a(n) - a(2^k - n) = k2^{k-1} - k(2^k - n) = k(n - 2^{k-1})$.

This proves (4.8) by induction since the initial values $a(0) = a(1) = 0$ and $a(2) = 1$ give $d_0 = 1$, $d_1 = 1$, and $d_2 = -x^2$.

For example $a(3)$ is the number of 1's in $1, 10$, i.e. $a(3) = 2$. Thus

$$d_3 = \det\begin{pmatrix} 1 & x & 0 \\ x & 0 & x^2 \\ 0 & x^2 & 0 \end{pmatrix} = -x^{4+0} = x^4.$$

**Example 4.8**

Let $c(2^k - 1) = x^{2^k-1}$ and $c(n) = 0$ else. Then

$$\mathbf{d}_n = \det\left(c(i+j)\right)_{i,j=0}^{n-1} = (-1)^{\binom{n}{2}} x^{2\binom{n}{2}}. \tag{4.9}$$

For example

$$\mathbf{d}_4 = \det\begin{pmatrix} 1 & x & 0 & x^3 \\ x & 0 & x^3 & 0 \\ 0 & x^3 & 0 & 0 \\ x^3 & 0 & 0 & 0 \end{pmatrix} = x^{12} = x^{2\binom{4}{2}}.$$

For the proof observe that $\mathbf{d}_n = (-1)^n x^{(2^k-1)(2n-2^k)} \mathbf{d}_{2^k-n}$ for $1 < 2^{k-1} < n \leq 2^k$.

It suffices to verify that $2\binom{2^k-n}{2} + (2^k-1)(2n-2^k) = 2\binom{n}{2}$.

**5. The matrices $\left(a_{i+j+1}\right)_{i,j=0}^{n-1}$.**

**Theorem 5.1**

Let $D(n) = \det\left(a_{i+j+1}\right)_{i,j=0}^{n-1}$, $\gamma(n) = 2^{\lfloor \log_2(n) \rfloor} - 1$ and $\delta(n) = 2n - \gamma(n)$.

Then

$$D(n) = x_{\gamma(n)}^{\delta(n)} (-1)^n D(\gamma(n) - n). \tag{5.1}$$

**Proof**

For given $n > 0$ choose $k$ such that $2^{k-1} < n+1 \leq 2^k$. Then $k-1 < \log_2(n+1) \leq k$ and $\gamma(n) = 2^k - 1$.

Let $\pi$ be a 1-nimble permutation. Then as above we see that $\pi$ induces an order reversing permutation on the interval $[2^k - 1 - n, n-1]$. Here we have $i + \pi(i) = 2^k - 2$.

Since there are $\delta(n) = 2n - \gamma(n)$ elements in the interval $[2^k - 1 - n, n-1]$ we get (5.1) by induction.

**Corollary 5.2**

The sequence $D(n)$ satisfies

$$D(2^k + n) = (-1)^n x_{2^{k+1}-1}^{2n+1} D(2^k - 1 - n) \tag{5.2}$$

for $0 \leq n < 2^k$.

Let us compute the first values with this recursion:

$D(0) = 1 \quad D(1) = x_1$
$D(3) = -x_3^3 \quad D(2) = x_1 x_3$

$D(0) = 1 \quad D(1) = x_1 \quad D(2) = x_1 x_3 \quad D(3) = -x_3^3$
$D(7) = -x_7^7 \quad D(6) = x_1 x_7^5 \quad D(5) = -x_1 x_3 x_7^3 \quad D(4) = -x_3^3 x_7$

For example for $n = 2^{k+1} - 1$ we have $\gamma(n) = n$ and $\delta(n) = n$.

This implies for $k > 1$

$$D(2^k - 1) = -x_{2^k - 1}^{2^k - 1}. \tag{5.3}$$

**Lemma 5.3**

For $0 \leq n < 2^{k-1}$ we get

$$D(2^k + n) = -x_{2^{k+1}-1}^{2n+1} x_{2^k-1}^{2^k-1-2n} D(n). \tag{5.4}$$

**Proof**

By (5.2) we get

$$D(2^k + n) = (-1)^n x_{2^{k+1}-1}^{2n+1} D(2^k - 1 - n) = (-1)^n x_{2^{k+1}-1}^{2n+1} D(2^{k-1} + 2^{k-1} - 1 - n)$$

Again by (5.2) we have $D(2^{k-1} + 2^{k-1} - 1 - n) = (-1)^{n+1} x_{2^k-1}^{2^k-1-2n} D(n)$.

Thus $D(2^k + n) = -x_{2^{k+1}-1}^{2n+1} x_{2^k-1}^{2^k-1-2n} D(n)$.

**Lemma 5.4**

For $2^{k-1} \leq n < 2^k$ we get

$$D(2^k + n) = x_{2^{k+1}-1}^{2n+1} x_{2^k-1}^{2^k-1-2n} D(n). \tag{5.5}$$

**Proof**

$D(2^k + 2^{k-1} + i) = (-1)^i x_{2^{k+1}-1}^{2^k+2i+1} D(2^k - 1 - 2^{k-1} - i) = (-1)^i x_{2^{k+1}-1}^{2^k+2i+1} D(2^{k-1} - 1 - i)$
$= (-1)^i x_{2^{k+1}-1}^{2^k+2i+1} (-1)^i x_{2^k-1}^{-2i-1} D(2^{k-1} + i)$

which is equivalent with (5.5).

Let us write

$$D(n) = (-1)^{\sum_i \binom{\mu_i(n)}{2}} \prod_{i \geq 1} x_{2^i-1}^{\mu_i(n)}. \tag{5.6}$$

Then $\mu_k(n)$ only depends on the residue class modulo $2^{k+2}$.

**Theorem 5.5**

Let $0 \leq i \leq 2^{k-1} - 1$. Then

$$\begin{aligned}\mu_k(i) &= 0, \\ \mu_k(2^{k-1} + i) &= 2i + 1, \\ \mu_k(2^k + i) &= 2^k - 2i - 1, \\ \mu_k(2^k + 2^{k-1} + i) &= 0.\end{aligned} \qquad (5.7)$$

For $n \geq 2^{k+1}$ we have $\mu_k(n) = \mu_k(n \bmod 2^{k+1})$.

**Proof**

Formula (5.7) follows from (5.1).

By (5.4) we have $\mu_k(2^R + n) = \mu_k(n)$ for $R \geq k+1$ which implies $\mu_k(n) = \mu_k(n \bmod 2^{k+1})$.

Thus for $n < 2^k$ we have $\mu_k(i) =$

Let us for example compute $D(11)$.

$11 \equiv 3 = 2 + 1 \bmod 2^2$ implies $\mu_1(11) = \max(2 - 1 - 2, 0) = 0$,

$11 \equiv 3 = 2 + 1 \bmod 2^3$ implies $\mu_2(11) = \max(2 + 1, 0) = 3$,

$11 \equiv 11 = 2^3 + 3 \bmod 2^4$ implies $\mu_3(11) = \max(2^3 - 1 - 6, 0) = 1$,

$11 \equiv 11 = 8 + 3 \bmod 2^5$ implies $\mu_4(11) = \max(6 + 1, 0) = 7$.

Therefore we get $D(11) = (-1)^{\binom{3}{2}+\binom{1}{2}+\binom{7}{2}} x_3^3 x_7 x_{15}^7 = x_3^3 x_7 x_{15}^7$.

Let us now determine the numbers

$$T_n = \frac{D(n)D(n+2)}{D(n+1)^2}. \qquad (5.8)$$

From

$$D(2^k + 2^{k-1} - 2) = -x_{2^{k+1}-1}^{2^k-3} x_{2^k-1}^3 D(2^{k-1} - 2),$$
$$D(2^k + 2^{k-1} - 1) = -x_{2^{k+1}-1}^{2^k-1} x_{2^k-1}^1 D(2^{k-1} - 1),$$
$$D(2^k + 2^{k-1}) = x_{2^{k+1}-1}^{2^k+1} x_{2^k-1}^{-1} D(2^{k-1}),$$
$$D(2^k + 2^{k-1} + 1) = x_{2^{k+1}-1}^{2^k+3} x_{2^k-1}^{-3} D(2^k + 1)$$

we get

$T_{2^k + 2^{k-1} - 2} = -T_{2^{k-1} - 2}$ and $T_{2^k + 2^{k-1} - 1} = -T_{2^{k-1} - 1}$.

For $m > 0$ and $0 \leq j \leq 3$ we get by (5.4)

$$D\left(2^{k+m} + 2^{k-1} - 2 + j\right) = -x_{2^{k+m+1}-1}^{2^k-3+2j} x_{2^{k+m}-1}^{2^{k+m}+3-2^k-2j} D\left(2^{k-1} - 2 + j\right)$$

which implies $T_{2^{k+m}+2^{k-1}-2} = T_{2^{k-1}-2}$ and $T_{2^{k+m}+2^{k-1}-1} = T_{2^{k-1}-1}$.

The same argument using (5.4) gives $T_{2^{k+m}+2^k+2^{k-1}-2+j} = T_{2^k+2^{k-1}-2+j} = -T_{2^{k-1}-2+j}$ for $m > 1$ and $j = 0, 1$.

There remains $T_{2^{k+2}+2^{k+1}+2^k-2+j}$. By (5.5) we get $T_{2^{k+2}+2^{k+1}+2^k-2+j} = T_{2^{k+1}+2^k-2+j} = -T_{2^k-2+j}$.

This gives

**Theorem 5.6**

*The numbers $T_{2^{k+1}n+2^k-2+j}$ satisfy*

$$T_{2^{k+1}n+2^k-2+j} = (-1)^n T_{2^k-2+j} \qquad (5.9)$$

*for $j = 0, 1$.*

The first terms of the sequence $(T_n)_{n \geq 0}$ are

$$\frac{x_3}{x_1}, \quad -\frac{x_3}{x_1}, \quad -\frac{x_1 x_7}{x_3^2}, \quad \frac{x_1 x_7}{x_3^2}, \quad -\frac{x_3}{x_1}, \quad \frac{x_3}{x_1}, \quad -\frac{x_1 x_{15}}{x_7^2}, \quad \frac{x_1 x_{15}}{x_7^2}, \ldots.$$

**Theorem 5.7**

*The numbers $T_n$, $n \geq 0$, satisfy*

$$\begin{aligned}
T_{2n+1} &= -T_{2n}, \\
T_{4n} &= (-1)^n \frac{x_3}{x_1} \\
T_{2^{k+1}n+2^k-1} &= (-1)^n \frac{x_1 x_{2^{k+1}-1}}{x_{2^k-1}^2} \quad \text{for } k \geq 2.
\end{aligned} \qquad (5.10)$$

**Proof**

For $2^{k-1} \leq n < 2^k - 2$ we have

$$D(n) = (-1)^n x_{2^k-1}^{2n+1-2^k} D\left(2^k - n - 1\right),$$

$$D(n+1) = (-1)^{n+1} x_{2^k-1}^{2n+3-2^k} D\left(2^k - n - 2\right),$$

$$D(n+2) = (-1)^n x_{2^k-1}^{2n+5-2^k} D\left(2^k - n - 3\right).$$

This implies

$$\frac{D(n)D(n+2)}{D(n+1)^2} = \frac{x_{2^k-1}^{2n+1-2^k} D(2^k-n-1) x_{2^k-1}^{2n+5-2^k} D(2^k-n-3)}{x_{2^k-1}^{2n+3-2^k} D(2^k-n-2) x_{2^k-1}^{2n+3-2^k} D(2^k-n-2)} = \frac{D(2^k-n-1)D(2^k-n-3)}{D(2^k-n-2)^2}.$$

Therefore we have

$$T_n = T_{2^k-3-n} \qquad (5.11)$$

for $2^{k-1} \le n \le 2^k - 3$.

It remains to prove

$$T_{2^k-2} = -\frac{x_1 x_{2^{k+1}-1}}{x_{2^k-1}^2},$$

$$T_{2^k-1} = \frac{x_1 x_{2^{k+1}-1}}{x_{2^k-1}^2} \qquad (5.12)$$

for $k \ge 2$.

By (5.4) we have

$$D(2^k+1) = -x_1 x_{2^{k+1}-1}^3 x_{2^k-1}^{2^k-3},$$
$$D(2^k) = -x_{2^{k+1}-1} x_{2^k-1}^{2^k-1}.$$

By (5.3) $D(2^k-1) = -x_{2^k-1}^{2^k-1}$ and by (5.2) $D(2^k-2) = x_1 x_{2^k-1}^{2^k-3}$.

This gives (5.12).

Let us prove that $T_{4n} = (-1)^n \frac{x_3}{x_1}$ and $T_{4n+1} = (-1)^{n+1} \frac{x_3}{x_1}$.

By induction we get using (5.11)

$$T_{4(2^k+j)} = T_{2^{k+3}-3-4(2^k+j)} = T_{4(2^k-j-1)+1} = (-1)^{2^k+j} \frac{x_3}{x_1}.$$

$$T_{4(2^k+j)+1} = T_{2^{k+3}-3-4(2^k+j)-1} = T_{4(2^k-j-1)} = (-1)^{2^k+j+1} \frac{x_3}{x_1}.$$

**Remarks**

Let us derive some connections between $D(n)$ and $d(n)$.

**Lemma 5.8**

$$\frac{d(n)d(n+1)}{D(n)^2} = (-1)^n x_0. \qquad (5.13)$$

**Proof**

For $0 < n < 2^k$ we have by (4.2)

$$d(2^k + n) = (-1)^n x_{2^{k+1}-1}^{2n} d(2^k - n),$$
$$d(2^k + n + 1) = (-1)^{n+1} x_{2^{k+1}-1}^{2n+2} d(2^k - n - 1)$$

and therefore

$$d(2^k + n) d(2^k + n + 1) = -x_{2^{k+1}-1}^{4n+2} d(2^k - n) d(2^k - n - 1).$$

By (5.2) we have

$$D(2^k + n)^2 = x_{2^{k+1}-1}^{4n+2} D(2^k - 1 - n)^2.$$

This implies

$$h(2^k + n) = \frac{d(2^k + n) d(2^k + n + 1)}{D(2^k + n)^2} = \frac{-d(2^k - n) d(2^k - n - 1)}{D(2^k - 1 - n)^2} = -h(2^k - n - 1)$$

for $0 < n < 2^k$. Further we have

$$h(2^k) = \frac{d(2^k) d(2^k + 1)}{D(2^k)^2} = \frac{d(2^k)(-1)x_{2^{k+1}-1}^2 d(2^k - 1)}{x_{2^{k+1}-1}^2 D(2^k - 1)^2} = -h(2^k - 1).$$

Therefore we get

$$h(2^k + n) = -h(2^k - n - 1)$$

for $0 \leq n < 2^k$ and $k \geq 1$. This gives by induction
$$h(2^k + n) = -h(2^k - 1 - n) = (-1)^{2^k - n} x_0 = (-1)^{2^k + n} x_0.$$

**Example 5.9**

*Let $b_{2^k - 1} = x^k$ and $b_n = 0$ else. Then*

$$D_n = \det(b_{i+j+1})_{i,j=0}^{n-1} = (-1)^{\delta(n)} x^{a(n) + a(n+1)}, \qquad (5.14)$$

*where $a(n)$ denotes the total number of 1's in the binary expansions of the first $n-1$ positive integers.*

This follow immediately from (4.8) and (5.13).

For the numbers $T_n$ we get

$$T_n = (-1)^{\tau_n} x^{s_2(n+2) - s_2(n)}, \qquad (5.15)$$

where $s_2(n)$ denotes the sum of digits of the binary expansion of $n$.

For

$$T_n = \frac{D(n)D(n+2)}{D(n+1)^2} = (-1)^{\tau_n} \frac{x^{a(n)+a(n+1)} x^{a(n+2)+a(n+3)}}{x^{2a(n+1)+2a(n+2)}} = (-1)^{\tau_n} x^{a(n+3)-a(n+2)-(a(n+1)-a(n))}$$
$$= (-1)^{\tau_n} x^{s_2(n+2)-s_2(n)}.$$

The first terms of $(T_n)_{n \geq 0}$ are $x, -x, -1, 1, -x, x, -\frac{1}{x}, \frac{1}{x}, x, -x, 1, -1, -x, x, -\frac{1}{x^2}, \frac{1}{x^2}, x, -x, \cdots$.

By (5.10) we have

$T_{2n+1} = -T_{2n},$

$T_{4n} = (-1)^n x$

$T_{2^{k+1}n+2^k-1} = (-1)^n \frac{x^{k+2}}{x^{2k}} = (-1)^n x^{2-k}$ for $k \geq 2$.

This is in accord with

$$s_2(2^{k+1}n + 2^k + 1) - s_2(2^{k+1}n + 2^{k-1} + 2^{k-2} + \cdots + 1) = 2 - k.$$

**Example 5.10**

Let $c(2^k - 1) = x^{2^k - 1}$ and $c(n) = 0$ else. Then

$$\mathbf{D}_n = \det(c(i+j+1))_{i,j=0}^{n-1} = (-1)^{\tau_n} x^{n^2}. \tag{5.16}$$

**Proof**

We know that $\mathbf{d}_n = (-1)^{\binom{n}{2}} x^{2\binom{n}{2}}$. Since $\binom{n}{2} + \binom{n+1}{2} = n^2$ we get (5.16).

In this case we get $T_n = (-1)^{\tau_n} x^2$.

**Theorem 5.11**

If we choose $g(1) = 1$, and $g(2^k - 1) = (-1)^k$ for $k > 1$ and $g(n) = 0$ else then

$$\det(g(i+j+1))_{i,j=0}^{n-1} = r(n) \tag{5.17}$$

where $r(n)$ denotes the Golay-Rudin-Shapiro sequence, which is defined by

$$\begin{aligned} r(2n) &= r(n), \\ r(2n+1) &= (-1)^n r(n), \\ r(0) &= 1. \end{aligned} \tag{5.18}$$

**Proof**

Some information about the Golay-Rudin-Shapiro sequence can be found in OEIS [8] A020985. As has been observed in [2] the Golay-Rudin-Shapiro sequence counts the number of pairs 11 in the binary expansion of $n$ modulo 2:

$$r(n) = (-1)^{\varepsilon_0 \varepsilon_1 + \cdots + \varepsilon_{k-1} \varepsilon_k} \text{ if } n = [\varepsilon_k \varepsilon_{k-1} \cdots \varepsilon_0]_2. \qquad (5.19)$$

Thus $r(0) = (-1)^0 = 1$.

If $n = [\varepsilon_k \varepsilon_{k-1} \cdots \varepsilon_0]_2$ then $2n = [\varepsilon_k \varepsilon_{k-1} \cdots \varepsilon_0 0]_2$ and $2n+1 = [\varepsilon_k \varepsilon_{k-1} \cdots \varepsilon_0 1]_2$

which implies $r(2n) = r(n)$ and $r(2n+1) = r(n)(-1)^n$ because $\varepsilon_0 1 = 11$ if $n$ is odd.

The Golay-Rudin-Shapiro sequence can also be characterized by the recursion

$$r(2^k + n) = r(n) \text{ for } 0 \le n < 2^{k-1},$$
$$r(2^k + n) = -r(n) \text{ for } 2^{k-1} \le n < 2^k \qquad (5.20)$$
$$\text{for } k \ge 2 \text{ and } r(0) = r(1) = 1.$$

If $n < 2^{k-1}$ then $2^k + n = [10 \cdots]_2$ and thus $r(2^k + n) = r(n)$.

If $2^{k-1} \le n < 2^k$ then $2^k + n = [11 \cdots]_2$ and thus $r(2^k + n) = -r(n)$.

The proof of (5.17) now follows from (5.4) and (5.5).

From (2.10) we get

**Corollary 5.12**

$$\sum_{k \ge 0} (-1)^k z^{2^k - 1} = \cfrac{1}{1 + \cfrac{r(0)r(2)z}{1 + \cfrac{r(1)r(3)z}{1 + \cfrac{r(2)r(4)z}{1 + \ddots}}}}. \qquad (5.21)$$

## 6. Hankel determinants of shifted sequences $(a_{n+m})_{n \ge 0}$.

All results are very similar to the case $x_n = 1$. Therefore we only need to make slight alterations.

Let us state the first terms of $d(n,m) = \det(a_{i+j+m})_{i,j=0}^{n-1}$ for $m = 3$ and $m = 5$:

$$(d(n,3)) = \left(1, x_3, 0, 0, -x_3 x_7^3, x_7^5, 0, 0, -x_7^5 x_{15}^3, -x_3 x_7^3 x_{15}^5, 0, 0, -x_3 x_{15}^{11}, \cdots \right),$$

$$(d(n,5)) = \left(1, 0, 0, -x_7^3, 0, 0, 0, 0, -x_7^3 x_{15}^5, 0, 0, -x_{15}^{11}, \cdots \right).$$

As in Lemma 3.3 we see that

$$d(n,m) = \det\left(a_{i+j+m}\right)_{i,j=0}^{n-1} = 0 \qquad (6.1)$$

if $2^k - m < n < 2^k$ for some $k$.

**Lemma 6.1**

Let $m = 2r+1$ and $k$ be given. Then for $0 \leq j \leq 2^{k-1} - r - 1$

$$d\left(2^{k-1} - r + j, m\right) = (-1)^{\binom{2j+1}{2}} x_{2^k-1}^{2j+1} d\left(2^{k-1} - r - j - 1, m\right). \qquad (6.2)$$

For example $d(9,3) = d(8-1+2,3) = (-1)^{10} x_{15}^5 d(8-1-2-1,3) = x_{15}^5 d(4,3)$.

**Lemma 6.2**

Let $m = 2r+1$. Then for $a \leq 2^R - m$ for some $R$ we have

$$d\left(2^R + a, m\right) = (-1)^{\binom{2a+2r+1}{2}} x_{2^{R+1}-1}^{2a+m} d\left(2^R - a - m, m\right). \qquad (6.3)$$

if $a \leq 2^R - m$.

**Lemma 6.3**

Let $m = 2r$ and $k$ be given. Then for $0 \leq j \leq 2^{k-1} - r$

$$d\left(2^{k-1} - r + j, m\right) = (-1)^{\binom{2j}{2}} x_{2^k-1}^{2j} d\left(2^{k-1} - r - j, m\right). \qquad (6.4)$$

**Lemma 6.4**

Let $m = 2r$. Then for $a \leq 2^R - m$ for some $R$ we have

$$d\left(2^R + a, m\right) = (-1)^{\binom{2a+m}{2}} x_{2^{R+1}-1}^{2a+m} d\left(2^R - a - m, m\right). \qquad (6.5)$$

**Theorem 6.5**

$d(n,m) \neq 0$ if and only if $n \equiv 0, -m \bmod 2^{K+1}$ if $2^K < m \leq 2^{K+1}$.